\pdfoutput=1

\documentclass[12pt]{amsart} 
\usepackage{graphicx,color}
\usepackage[margin=2cm]{geometry}
\usepackage{amsthm,amsmath,amssymb,enumerate,comment}
\usepackage{mathtools}

\usepackage{cleveref}
\usepackage{autonum}
\usepackage[all]{xy}
\usepackage{here}
\usepackage{tikz}
\usetikzlibrary{calc}
\usetikzlibrary {arrows.meta}

\newcommand{\F}{\mathbb{F}}

\newcommand{\Z}{\mathbb{Z}}

\newcommand{\confi}[2]{C_{#1}({#2})}
\newcommand{\uconfi}[2]{B_{#1}({#2})}
\newcommand{\torus}[1]{T^{#1}}
\newcommand{\sphere}[1]{S^{#1}}
\newcommand{\rp}[1]{\mathbb{R}P^{#1}}
\newcommand{\sym}[1]{\Sigma_{#1}}

\newcommand{\cohomology}[2]{H^{#1}({#2})}
\newcommand{\ccohomology}[3]{H^{#1}({#2};{#3})}

\newcommand{\lcohomology}[3]{\mathcal{H}^{#1}({#2};{#3})}

\newcommand{\abhtpygrp}[2]{\pi_{#1}({#2})}

\newcommand{\map}[3]{{#1}\colon{#2}\rightarrow {#3}}
\newcommand{\namelessmap}[2]{{#1}\rightarrow {#2}}

\newcommand{\iset}[2]{\{{#1}|{#2}\}}

\newcommand{\Hom}[3]{\text{Hom}_{#1}({#2},{#3})}

\newtheoremstyle{mystyle}
    {3pt}
    {3pt}
    {\itshape}
    {}
    {\bfseries}
    {}
    {10pt}
    {\thmname{#1} \thmnumber{#2}\thmnote{(#3)}}

\theoremstyle{mystyle}

\newtheorem{theorem}{Theorem}[section]
\newtheorem{proposition}{Proposition}[section]
\newtheorem{lemma}{Lemma}[section]
\newtheorem{corollary}{Corollary}[section]

\numberwithin{equation}{section}

\makeatletter
\renewenvironment{proof}[1][\proofname]{\par
  \pushQED{\qed}%
  \normalfont \topsep6\p@\@plus6\p@\relax
  \trivlist
  \item\relax
  {#1\@addpunct{.}}\hspace\labelsep\ignorespaces
}{%
  \popQED\endtrivlist\@endpefalse
}
\makeatother

\renewcommand{\proofname}{\textbf{Proof}}

\crefname{section}{Chapter}{Chps.}
\crefname{corollary}{Corollary}{Cors.}
\crefname{equation}{eq.}{eqs.}
\crefname{align}{eq.}{eqs.}
\crefname{figure}{Fig.}{figs.}
\crefname{theorem}{Theorem}{Theorems}
\crefname{lemma}{Lemma}{Lems.}

\title{Mod 2 Cohomology of 2-configuration space of a closed surface and Stiefel--Whitney class}
\author{Tomoki TOKUDA}
\email{tokuda.tomoki.379@s.kyushu-u.ac.jp}
\address{Joint Graduate School of Mathematics for Innovation, Kyushu University, Fukuoka, 819-0395, Japan}

\begin{document}
\maketitle
\markboth{TOMOKI TOKUDA}{MOD 2 COHOMOLOGY OF 2-CONFIGURATION SPACE OF A CLOSED SURFACE}

\begin{abstract}
    In this paper, we compute the singular cohomology groups $\ccohomology{i}{\confi{2}{M}}{\F_2}$ of the ordered 2-configuration space $\confi{2}{M}$ as $\sym{2}$-representations. Using the result, we determine the mod 2 cohomology of the unordered 2-configuration space $\uconfi{2}{M}$ as a $\ccohomology{*}{\rp{\infty}}{\F_2}$-module. As a corollary of our computation, we see that the Stiefel--Whitney height of $\confi{2}{M}$ is $2$ or $3$ when $M$ is orientable or not, respectively.
\end{abstract}

\section{Introduction}
The ordered $n$-configuration space $\confi{n}{X}$ of a space $X$ is the subspace of the $n$-fold product space $X^{n}$ defined by 
\begin{align}
    \confi{n}{X}=\iset{(x_1,\dots ,x_n)\in X^{n}}{x_i\neq x_j \text{ if } i\neq j}.
\end{align}
The symmetric group $\sym{n}$ acts freely on $\confi{n}{X}$. The orbit space $\uconfi{n}{X}=\confi{n}{X}/\sym{n}$ is called the unordered configutation space. For a commutative ring $R$, the singular cohomology $\ccohomology{*}{\confi{n}{X}}{R}$ is a $\sym{n}$-representation. The main result of this paper is calculation of the $\sym{2}$-representation $\ccohomology{i}{\confi{2}{M}}{\F_2}$ for every closed surface $M$ (\cref{c_ori,c_unori}).

\vspace{3mm}
In \cite{tokuda2024mod2representationsymmetric}, the author determined the $\ccohomology{*}{\rp{\infty}}{\F_2}$-module structure of $\ccohomology{*}{\uconfi{2}{\torus{d}}}{\F_2}$ for $d=2,3$ where $\torus{d}$ is the $d$-torus. The computation of \cite{tokuda2024mod2representationsymmetric} uses the Serre spectral sequence associated to the homotopy fibration
\begin{align}
    \xymatrix{
    \confi{2}{\torus{d}} \ar[r] &\uconfi{2}{\torus{d}}\ar[r] &\rp{\infty}.
    }
\end{align}
This homotopy fibration is obtained by the Borel construction
\begin{align}
    \xymatrix{
    \confi{2}{\torus{d}}\ar[r] &\confi{2}{\torus{d}}\times_{\sym{2}}\sphere{\infty} \ar[r]& \rp{\infty}.
    }
\end{align}
and canonical homotopy equivalence $\namelessmap{\confi{2}{\torus{d}}\times_{\sym{2}}\sphere{\infty}}{\uconfi{2}{\torus{d}}}$. 

\vspace{5mm}
For every closed surface $M$, we may consider the similar homotopy fibration
\begin{align}\label{htpy_fib1}
    \xymatrix{
    \confi{2}{M}\ar[r] &\uconfi{2}{M}\ar[r] &\rp{\infty}.
    }
\end{align}
Let $E$ be the Serre spectral sequence associated to (\ref{htpy_fib1}). Then the $E_\infty$-term is isomorphic to $\ccohomology{*}{\uconfi{2}{M}}{\F_2}$ as a $\ccohomology{*}{\rp{\infty}}{\F_2}$-module. Thus computing $E$, we determined the $\ccohomology{*}{\rp{\infty}}{\F_2}$-module structure of $\ccohomology{i}{\uconfi{2}{M}}{\F_2}$ (\cref{b_ori,b_unori}). 

An advantage of this calculation method is that we can see the Stiefel--Whitney class of $\confi{2}{M}$. An explicit definition of the Stiefel--Whitney class of a free $\sym{2}$-space is introduced in \cite{CAT}, Defintion 8.21.

\vspace{10mm}
The point of this study is the computation of the cohomology groups $\ccohomology{i}{\confi{2}{M}}{\F_2}$ as the $\sym{2}$-representations. For the case of torus, we easily see that the ordered 2-configuration space $\confi{2}{\torus{d}}$ is homeomorphic to the product $\torus{d}\times (\torus{d}\setminus *)$ by the group structure of torus. However, most of closed surfaces do not equip a group structure, so $\confi{2}{M}\not\cong M\times (M\setminus *)$ in general. 

So we first compute $\ccohomology{i}{\confi{2}{M}}{\F_2}$ again using the Serre spectral sequence $E'$ associated to the Fadell--Neuwirth fibration
\begin{align}\label{htpy_fib2}
    \xymatrix{
    M\setminus * \ar[r]& \confi{2}{M}\ar[r]&M.
    }
\end{align}
where the projection $\namelessmap{\confi{2}{M}}{M}$ is the 1st coodinate projection, $(x,y)\mapsto x$. 

\vspace{5mm}
The second pages $E_2$ and $E_{2}^{\prime}$ of the spectral sequence $E$ and $E^{\prime}$ are given by 
\begin{align}
    E_{2}^{p,q}&=\ccohomology{p}{\rp{\infty}}{\lcohomology{q}{\confi{2}{M}}{\F_2}},\\
    E_{2}^{\prime s,t}&=\ccohomology{s}{M}{\lcohomology{t}{M\setminus *}{\F_2}}
\end{align}
where $\mathcal{H}$ denotes the local system associated to each fibration.

\vspace{5mm}
By the classification theorem of closed surfaces, we consider the cases $M$ is the sphere $\sphere{2}$, $g$-fold connected sum $T\#\cdots\# T$ of tori, or the $k$-fold connected sum $\rp{2}\#\cdots \#\rp{2}$. Here are the result of this paper.

\vspace{3mm}
\begin{theorem}\label{c_ori}
    When $M$ is an orientable closed surface with genus $g$, the $\sym{2}$-representations $\ccohomology{i}{\uconfi{2}{M}}{\F_2}$ decompose as follows:
    \begin{align}
        \ccohomology{i}{\confi{2}{M}}{\F_2}\cong \begin{dcases}
            \F_2 & \text{$i=0$,}\\
            \F_2[\sym{2}]^{\oplus 2g} & \text{$i=1$,}\\
            \F_{2}^{\oplus 2g+1} \oplus \F_2[\sym{2}]^{\oplus 2g^2+g} & \text{$i=2$,}\\
            \F_{2}^{\oplus 2g} & \text{$i=3$,}\\
            0 & \text{otherwise.}
        \end{dcases}
    \end{align}
\end{theorem}

\vspace{3mm}
Let $\alpha$ be the generator of the cohomology ring $\ccohomology{*}{\rp{\infty}}{\F_2}$. The induced map of the projection $\map{p}{\uconfi{2}{M}}{\rp{\infty}}$, in \ref{htpy_fib1}, maps $\alpha\in\ccohomology{1}{\rp{\infty}}{\F_2}$ to $p^*(\alpha)\ccohomology{1}{\uconfi{2}{M}}{\F_2}$. To simplify our notation, we denote this $p^*(\alpha)$ by $\alpha$. 

\vspace{2mm}
This $\alpha\in \ccohomology{1}{\uconfi{2}{M}}{\F_2}$ is called the 1st Stiefel--Whitney class of $\confi{2}{M}$ (see \cite{CAT}, Definition 8.21). 
\begin{theorem}\label{b_ori}
    When $M$ is an orientable closed surface with genus $g$, there is an isomorphism of $\ccohomology{*}{\rp{\infty}}{\F_2}$-modules
    \begin{align}
        \ccohomology{*}{\uconfi{2}{M}}{\F_2}\cong \F_2[\alpha]/(\alpha^3)\oplus \left(\bigoplus_{i=1}^{g}\F_2 x_i\right)\oplus\left(\bigoplus_{i=1}^{2g^2+g}\F_2 y_i \right)\oplus \left(\bigoplus_{i=1}^{2g}\F_2[\alpha]/(\alpha^2)z_i\right)
    \end{align}
    where $deg(\alpha)=deg (x_i)=1$, $deg (y_i)=2$ and $deg (z_i)=3$, and $\alpha$ is the 1st Stiefel--Whitney class of $\confi{2}{M}$.
\end{theorem}

\vspace{3mm}
\begin{theorem}\label{c_unori}
    When $M$ is the $k$-fold connected sum $\rp{2}\#\cdots \#\rp{2}$, the $\sym{2}$-representations $\ccohomology{i}{\uconfi{2}{M}}{\F_2}$ decompose as follows:
    \begin{align}
        \ccohomology{i}{\confi{2}{M}}{\F_2}\cong \begin{dcases}
            \F_2 & \text{$i=0$,}\\
            \F_2[\sym{2}]^{\oplus k} & \text{$i=1$,}\\
            \F_{2}^{\oplus k-1} \oplus \F_2[\sym{2}]^{\oplus \frac{1}{2}k(k+1)+1} & \text{$i=2$,}\\
            \F_{2}^{\oplus k} & \text{$i=3$,}\\
            0 & \text{otherwise.}
        \end{dcases}
    \end{align}
\end{theorem}

\vspace{3mm}
\begin{theorem}\label{b_unori}
    When $M$ is the $k$-fold connected sum $\rp{2}\#\cdots \#\rp{2}$, there is an isomorphism of $\ccohomology{*}{\rp{\infty}}{\F_2}$-modules 
    \begin{align}
        \ccohomology{*}{\uconfi{2}{M}}{\F_2}\cong \F_2[\alpha]/(\alpha^4)\oplus \left(\bigoplus_{i=1}^{k}\F_2 x_i\right)\oplus\left(\bigoplus_{i=1}^{\frac{1}{2}k(k+1)+1}\F_2 y_i\right) \oplus \left(\bigoplus_{i=1}^{k-1}\F_2[\alpha]/(\alpha^2)z_i\right)
    \end{align}
    where $deg(\alpha)=deg (x_i)=1$, $deg (y_i)=2$ and $deg (z_i)=3$, and $\alpha$ is the 1st Stiefel--Whitney class of $\confi{2}{M}$.
\end{theorem}

\vspace{3mm}
The Stiefel--Whitney height of $M$ is the maximal number $m$ such that $\alpha^m\neq 0$ where $\alpha$ is the 1st Stiefel--Whitney class of the $\sym{2}$-space $\confi{2}{M}$ (see \cite{CAT}, Definition 8.24). As the corollary of these theorems, we obtain the following.
\begin{corollary}
    The Stiefel--Whitney height $h(M)$ of a closed orientable surface $M$ is $2$ or $3$ when $M$ is orientable or not, respectively.
\end{corollary}
In \cite{kishimoto2023van}, Kishimoto and Matsushita generalized van Kampen-Flores theorem using Stiefel--Whitney height.

\vspace{5mm}
\section{Preparation for computation}

We start with the computation for $\lcohomology{t}{M\setminus *}{\F_2}$. Note that the base space $M$ of the fibration \ref{htpy_fib2} is a closed surface, so it is not simply connected unless $M$ is the sphere. However, in fact, the action of $\abhtpygrp{1}{M}$ on $\lcohomology{t}{M\setminus *}{\F_2}$ is always trivial. 
\begin{lemma}
    For any closed surface $M$, the action of $\abhtpygrp{1}{M}$ on $\lcohomology{t}{M\setminus *}{\F_2}$ is trivial for all $t$. 
\end{lemma}

\begin{proof}
This depends on the fact that the inclusion $\namelessmap{M\setminus *}{M}$ induces a surjection $\namelessmap{\ccohomology{t}{M}{\F_2}}{\ccohomology{t}{M\setminus *}{\F_2}}$ for any integer $t$, which is obvious by cellular cohomology. Actually, this implies that the map of fiber bundles
\begin{align}\label{fiber_map}
    \xymatrix@C=42pt{
    M\setminus * \ar[r] \ar_{inclusion}[d] & \confi{2}{M} \ar^{pr_1}[r] \ar^{i=inclusion}[d] & M \ar@{=}[d]\\
    M \ar[r] & M\times M \ar_{pr_1}[r] & M
    }
\end{align}
induces a surjection $\namelessmap{\lcohomology{t}{M}{\F_2}}{\lcohomology{t}{M\setminus *}{\F_2}}$ of $\F_2[\abhtpygrp{1}{M}]$-modules for each $t$. But the action of $\abhtpygrp{1}{M}$ on $\lcohomology{t}{M}{\F_2}$ is trivial since the bundle $\namelessmap{M\times M}{M}$ is trivial. 
\end{proof}

So we may regard $\lcohomology{t}{M\setminus *}{\F_2}$ just as the trivial coefficient cohomology group $\ccohomology{t}{M\setminus *}{\F_2}$. This makes $E^{\prime}$ easy to compute.

\vspace{7mm}
Computing $E^{\prime}$, we see the cohomology groups $\ccohomology{q}{\confi{2}{M}}{\F_2}$. But what we really want is the local system $\lcohomology{q}{\confi{2}{M}}{\F_2}$. So we want to see the $\F_2[\abhtpygrp{1}{\rp{\infty}}]$-module structure of $\ccohomology{q}{\confi{2}{M}}{\F_2}$.

\vspace{5mm}
Here we turn our attention to the homomorphism $\map{i^*}{\ccohomology{q}{M\times M}{\F_2}}{\ccohomology{q}{\confi{2}{M}}{\F_2}}$ induced by the inclusion $\map{i}{\confi{2}{M}}{M\times M}$. Since this inclusion is an equivariant map with respect to the action of $\abhtpygrp{1}{\rp{\infty}}\cong \sym{2}$ on $\confi{2}{M}$ and $M\times M$. Thus it induces homomorphisms $\map{i^*}{\ccohomology{q}{M\times M}{\F_2}}{\ccohomology{q}{\confi{2}{M}}{\F_2}}$ of $\F_2[\sym{2}]$-modules. Applying the Leray--Hirsch theorem to the two fiber bundles in \ref{fiber_map}, we see that these $i^*$ are surjective. So we can see the $\F_2[\sym{2}]$-module structure of $\ccohomology{q}{\confi{2}{M}}{\F_2}$ as a quatient of $\ccohomology{q}{M\times M}{\F_2}$. 

\vspace{5mm}
Let $\Delta(M)$ denote the diagonal set in $M\times M$. Then $\confi{2}{M}$ is the complement $M\times M \setminus \Delta(M)$. Since $i^*$ is surjective, the long exact sequence of mod 2 cohomology groups associated to the pair $(M\times M,\confi{2}{M})$ splits into the short exact sequences
\begin{align}
    \xymatrix{
    0&\cohomology{q}{\confi{2}{M}}\ar[l]&\cohomology{q}{M\times M}\ar_{i^*}[l]&\cohomology{q}{M\times M, M\times M\setminus \Delta(M)}\ar[l]&0\ar[l]
    }
\end{align}
Further, by the Thom isomorphism theorem, we have
\begin{align}
    \cohomology{q}{M\times M, M\times M\setminus \Delta(M)}\cong \cohomology{q-2}{M}.
\end{align}
In particular, when $q=0$ or $q=1$, $\cohomology{q}{M\times M, M\times M\setminus \Delta(M)}=0$ and $i^*$ is an isomorphism. And $\cohomology{2}{M\times M, M\times M\setminus \Delta(M)}\cong \F_2$ and its generator is taken to the diagonal class $u_0\in \ccohomology{2}{M\times M}{\F_2}$ of $M$. Then we have the following proposition.
\begin{proposition}
    There is an isomorphism of $\F_2[\sym{2}]$-modules
    \begin{align}\label{quotient}
        \cohomology{q}{\confi{2}{M}}\cong \cohomology{q}{M\times M}/{\ker i^*}.
    \end{align}
\end{proposition}
Based on the above considerations, we will move on to the actual computation.

\vspace{5mm}
\section{Computation}
Let us compute the spactral sequences $E^{\prime}$ and $E$.  

\vspace{5mm}
\noindent
\textbf{Case 1}: $M$ is the sphere $\sphere{2}$.

In this case, $\uconfi{2}{M}$ is homotopy equivalent to the real projective plane $\rp{2}$ since there is a $\Z_2$-equivariant homotopy equivalence $\namelessmap{\sphere{2}}{\confi{2}{\sphere{2}}}$, $x\mapsto (x,-x)$. Indeed, this is a fiber homotopy equivalence between two double covers $\namelessmap{\sphere{2}}{\rp{2}}$, $\namelessmap{\confi{2}{\sphere{2}}}{\uconfi{2}{\sphere{2}}}$, so the induced map $\namelessmap{\rp{2}}{\uconfi{2}{\sphere{2}}}$ is a homotopy equivalence by the Whitehead theorem. Therefore, we do not need to calculate the Serre spectral sequences to see that
\begin{align}
    \ccohomology{*}{\uconfi{2}{\sphere{2}}}{\F_2}\cong \F_2[\alpha].
\end{align}

\vspace{5mm}
\noindent
\textbf{Case 2}: $M$ is the $g$-fold connected sum $T\#\cdots \# T$.

First we note that the integral cohomology algebra of $M=T\#\cdots \# T$ is the quotient of the polynomial algebra $\Z[a_1, \dots , a_g, b_1, \dots ,b_g]$ by the relations 
\begin{align}
    &a_i b_i= a_j b_j=u,\\
    &a_{i}^{2} = b_{i}^{2} = 0, \\
    &a_i a_j = b_i b_j = a_i b_j = 0 \qquad \text{if }i\neq j.
\end{align} Geometrically, each $a_i$ and $b_i$ correspond to the dual 1-cells of the meridian and longitude in the $i$-th torus. Thus the mod $2$ cohomology is similarly the quotient of $\F_2 [a_1, \dots , a_g, b_1, \dots ,b_g]$ by the same relations. 

\vspace{5mm}
Let us compute the spectral sequences $E'$ and $E$. Now $M\setminus *$ is $2g$-fold wedge sum of the circle. So the second page of $E'$ is given by 
\begin{align}\label{iso_of_bigraded}
    E^{\prime s,t}_{2}&\cong \ccohomology{s}{M}{\ccohomology{t}{M\setminus *}{\F_2}}\\
    &\cong \Hom{\Z}{\cohomology{s}{M}}{\ccohomology{t}{M\setminus *}{\F_2}}\\
    &\cong \cohomology{s}{M}\underset{\Z}{\otimes}\ccohomology{t}{M\setminus *}{\F_2}.
\end{align}

\begin{figure}[H]
    \centering
    \begin{tikzpicture}[xscale=1.5,yscale=1.5]
        \coordinate(O)at(0,0);
        \coordinate(XS)at(-0.5,0);
        \coordinate(XL)at(3.5,0);
        \coordinate(YS)at(0,-0.5);
        \coordinate(YL)at(0,2.5);
        \draw[semithick,->,>=stealth](XS)--(XL)node[right]{$s$};
        \draw[semithick,->,>=stealth](YS)--(YL)node[left]{$t$};
        \coordinate(P)at(1,0);
        \coordinate(Q)at(0,1);

        \foreach\k in{0,1,2,3}\foreach\l in{0,1,2}\fill($(O)+\k*(P)+\l*(Q)$)circle(0.05);

        \coordinate[label=above right:$1$] (P00) at (0,0);
        \coordinate[label=above right:$2g$, label=below :$a_i \otimes 1$] (P10) at (1,0);
        \coordinate[label= below: $b_i \otimes 1$] (Q10) at (1,-0.3);
        \coordinate[label=above right:$1$] (P20) at (2,0);
        \coordinate[label=above right:$0$] (P30) at (3,0);

        \coordinate[label=above right:$2g$, label=left: $1\otimes b_i$] (P01) at (0,1);
        \coordinate[label= left: $1 \otimes a_i\text{,}$] (Q10) at (-0.8,1);
        \coordinate[label=above right:$4g^2$] (P11) at (1,1);
        \coordinate[label=above right:$2g$] (P21) at (2,1);
        \coordinate[label=above right:$0$] (P31) at (3,1);

        \foreach\k in{0,1,2,3}
        \coordinate[label=above right:$0$](X\k)at($\k*(P)+(0,2)$);

        \draw [-{Stealth[length=3mm]}] (0.1,0.95) -- (1.9,0.1);
        \coordinate[label=above right:$d_2$](d)at(1.2,0.4);

    \end{tikzpicture}
    \caption{$E^{\prime}_{2}$-term}
    
\end{figure}
Note that \ref{iso_of_bigraded} is an isomorphism of bigraded algebra $E^{\prime}\cong \cohomology{*}{M}\otimes \ccohomology{*}{M\setminus *}{\F_2}$. So the Leibniz rule of the derivation implies that the $d_2$ in the figure above is zero. Hence $E^{\prime }$ collapses at $E^{\prime}_2$-term. 

The generators $a_i\otimes 1, 1\otimes a_i, b_i\otimes 1, 1\otimes b_i$ correspond to certain cross products in $\ccohomology{1}{M\times M}{\F_2}$, and they are taken to such generators of $\ccohomology{1}{\confi{2}{M}}{\F_2}$ by the inclusion $\namelessmap{\confi{2}{M}}{M\times M}$. Since the first cohomology $\ccohomology{1}{\confi{2}{M}}{\F_2}$ is isomorphic to $\ccohomology{1}{M\times M }{\F_2}$, we have an isomorphism of $\F_2[\sym{2}]$-modules
\begin{align}
    \ccohomology{1}{\confi{2}{M}}{\F_2}\cong \F_2[\sym{2}]^{\oplus 2g}
\end{align}
where the generators are $a_i\times 1, b_i \times 1,i=1,2,\dots, g$. 

By the ring structure of $\ccohomology{*}{M}{\F_2}$, the diagonal class $u_0\in \ccohomology{2}{M\times M,\confi{2}{M}}{\F_2}$ is written as $u_0=u\times 1+1\times u+\sum_{i} a_i\times b_i+\sum_{i}b_i\times a_i$. 

Let $v=u\times 1 + \sum_{i} a_i\times b_i$, then $u_0=v + \sigma v$ where $\sigma$ denotes the nontrivial element of $\sym{2}$. Since the second cohomology $\ccohomology{2}{M\times M}{\F_2}$ is the free module generated by $u\times 1, 1\times u, a_i\times a_j, a_i\times b_j, b_i\times a_j, b_i\times b_j, 1\le i,j\le g$, then \cref{quotient} implies 
\begin{align}
    \ccohomology{2}{\confi{2}{M}}{\F_2}\cong \F_2^{\oplus 2g+1}\oplus \F_2[\sym{2}]^{\oplus 2g^2+g}
\end{align}
where the $\F_2$-part is generated by $v,a_i\times a_i, b_i\times b_i, i=1,2,\dots, g$, and the $\F_2[\sym{2}]$-part is generated by $a_i\times a_j, b_i\times b_j$ for $1\le i<j \le g$ and $a_i\times b_j$ for $1\le i,j,\le g$.

\vspace{3mm}
The third cohomology $\ccohomology{3}{\confi{2}{M}}{\F_2}$ is the quotient of $\ccohomology{3}{M\times M}{\F_2}$ by the submodule generated by $(a_i\times 1)u_0, (1\times a_i)u_0, (b_i\times 1)u_0, (1\times b_i)u_0, i=1,2,\dots, g$. So, in $\ccohomology{3}{\confi{2}{M}}{\F_2}$, we have
\begin{align}
    u\times a_i + a_i\times u =&0,\\
    u\times b_i + b_i\times u =&0
\end{align}
for each $i=1,2,\dots, g$, and the generators do not have any other relation. Hence, 
\begin{align}
    \ccohomology{3}{\confi{2}{M}}{\F_2}\cong \F_{2}^{\oplus 2g}
\end{align}
as a $\F_2[\sym{2}]$-module where the generator is $u\times a_i, u\times b_i, i=1,2,\dots , g$. The fourth cohomology of $\confi{2}{M}$ is trivial since $\confi{2}{M}$ is a noncompact 4-manifold. This completes the proof of \cref{c_ori}.

\vspace{5mm}
Here we note that $\ccohomology{p}{\rp{\infty}}{\F_2[\sym{2}]}$ is isomorphic to $\ccohomology{p}{\sphere{\infty}}{\F_2}$. Now we can see the explicit form of $E$. 
\begin{align}\label{directsumdecomposition1}
    E_{2}^{p,q}&=\ccohomology{p}{\rp{\infty}}{\lcohomology{q}{\confi{2}{M}}{\F_2}}\\
    &=\begin{dcases}
        \ccohomology{p}{\rp{\infty}}{\F_2} & \text{$q=0$,}\\
        \ccohomology{p}{\rp{\infty}}{\F_2[\sym{2}]}^{\oplus 2g} & \text{$q=1$,}\\
        \ccohomology{p}{\rp{\infty}}{\F_2}^{\oplus 2g+1}\oplus \ccohomology{p}{\rp{\infty}}{\F_2[\sym{2}]}^{\oplus 2g^2+g} & \text{$q=2$,}\\
        \ccohomology{p}{\rp{\infty}}{\F_2}^{\oplus 2g} & \text{$q=3$,}
    \end{dcases}\\
    &=\begin{dcases}
        \ccohomology{p}{\rp{\infty}}{\F_2} & \text{$q=0$,}\\
        \ccohomology{p}{\sphere{\infty}}{\F_2}^{\oplus 2g} & \text{$q=1$,}\\
        \ccohomology{p}{\rp{\infty}}{\F_2}^{\oplus 2g+1}\oplus \ccohomology{p}{\sphere{\infty}}{\F_2}^{\oplus 2g^2+g} & \text{$q=2$,}\\
        \ccohomology{p}{\rp{\infty}}{\F_2}^{\oplus 2g} & \text{$q=3$.}
    \end{dcases}
\end{align}
Then $E_2$ is described as following. 
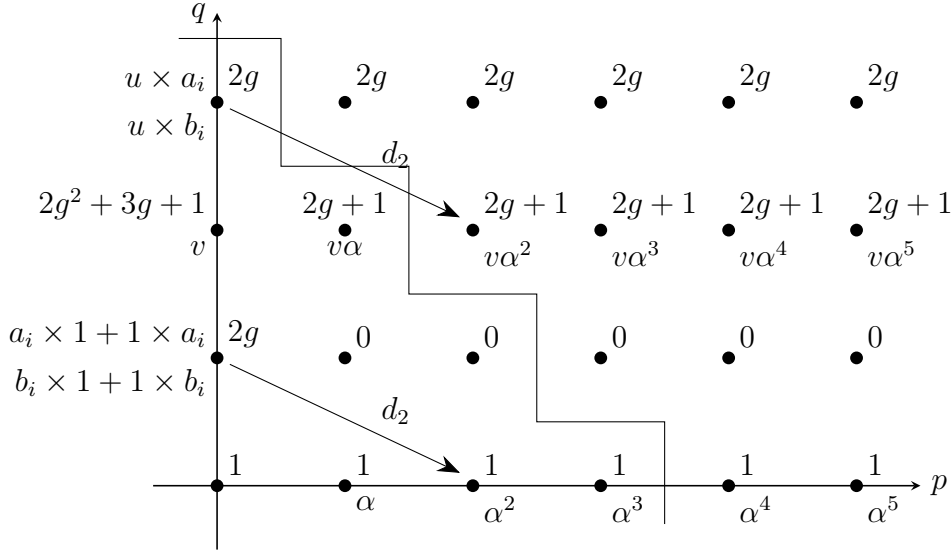
\begin{figure}[H]
    \centering
    \begin{tikzpicture}[xscale=1.7,yscale=1.7]
        \coordinate(O)at(0,0);
        \coordinate(XS)at(-0.5,0);
        \coordinate(XL)at(5.5,0);
        \coordinate(YS)at(0,-0.5);
        \coordinate(YL)at(0,3.7);
        \draw[semithick,->,>=stealth](XS)--(XL)node[right]{$p$};
        \draw[semithick,->,>=stealth](YS)--(YL)node[left]{$q$};
        \coordinate(P)at(1,0);
        \coordinate(Q)at(0,1);

        \foreach\k in{0,1,2,3,4,5}\foreach\l in{0,1,2,3}\fill($(O)+\k*(P)+\l*(Q)$)circle(0.05);
        \foreach\k in{0,1,2,3,4,5}
        \coordinate[label=above right:$1$](X\k)at($\k*(P)$);
        \coordinate[label=below right:$\alpha$](X)at(1,0);
        \foreach\k in{2,3,4,5}
        \coordinate[label=below right:$\alpha^{\k}$](X\k)at($\k*(P)$);
        \coordinate[label=above right:$2g$, label=above left:$a_i\times 1+1\times a_i$, label=below left:$b_i \times 1+1\times b_i$](Y)at(0,1);
        \foreach\k in{1,2,3,4,5}
        \coordinate[label=above right:$0$](Y\k)at($\k*(P)+(Q)$);
        \coordinate[label=above left:$2g^2+3g+1$, label=below left:$v$](Z)at(0,2);
        \coordinate[label=above:$2g+1$, label=below:$v \alpha$](Z1)at(1,2);
        \foreach\k in{2,3,4,5}\coordinate[label=above right:$2g+1$, label=below right:$v \alpha^{\k}$](Z\k)at($\k*(P)+(Z)$);
        \coordinate(W)at(0,3);
        \foreach\k in{0,1,2,3,4,5}\coordinate[label=above right:$2g$](W\k)at($\k*(P)+(W)$);
        \coordinate[label=above left:$u\times a_i$, label=below left:$u\times b_i$](G)at(0,3);

        \draw(-0.3,3.5)--(0.5,3.5)--(0.5,2.5)--(1.5,2.5)--(1.5,1.5)--(2.5,1.5)--(2.5,0.5)--(3.5,0.5)--(3.5,-0.3);

        \draw [-{Stealth[length=3mm]}] (0.1,0.95) -- (1.9,0.1);
        \coordinate[label=above right:$d_2$](d)at(1.2,0.4);
        \draw [-{Stealth[length=3mm]}] (0.1,2.95) -- (1.9,2.1);
        \coordinate[label=above right:$d_2$](d)at(1.2,2.4);

    \end{tikzpicture}
    \caption{$E_2$-term}
    
\end{figure}

Here we note that the $\alpha\in E_{2}^{1,0}$ is the Stiefel--Whitney class of $\confi{2}{M}$ at $E_{\infty}$.

\vspace{2mm}
Since $\uconfi{2}{M}$ is open 4-manifold, so the region outer the zigzag line in the figure must vanish at $E_{\infty}$. In particular, $E_{2}^{2,2}$ has to vanish at $E_{\infty}$-term. So the upper derivation in the figure above is injective. Here we remark that $u_0 \alpha^{2}$ is not hit by $d_2$. This is because the other generators, for example $(a_i\times a_i)\alpha^{2}$, cannot kill $\alpha^{5}$.

On the other hand, the lower $d_2$ in the figure must be zero because $d_2(a_i\times 1+1\times a_i)\neq 0$ implies $d_2((a_i\times 1+1\times a_i)\alpha)\neq 0$. Not only for elements $a_i\times 1+1\times a_i$ or $b_i\times 1+1\times b_i$, all elements in the $\ccohomology{*}{\sphere{\infty}}{\F_2}$-parts in the decomposition (\ref{directsumdecomposition1}) are taken to zero by the derivation $d_r$ for any $r\ge 2$. 

Then we obtain the $E_3$ as the following figure.

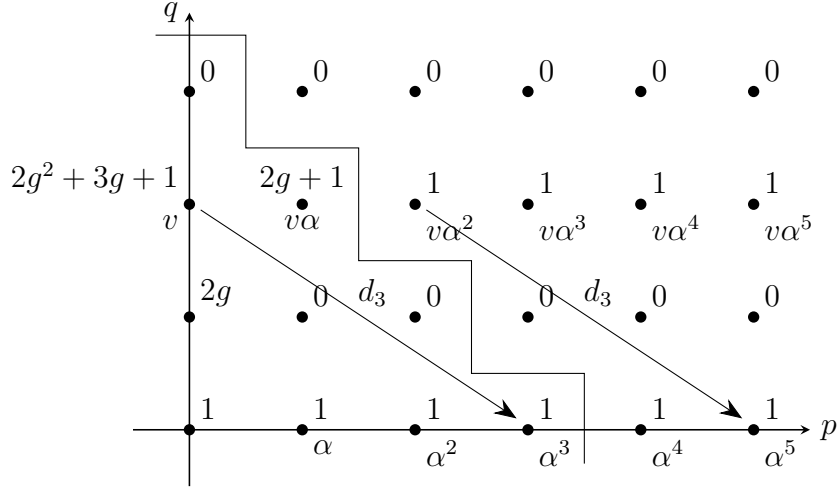
\begin{figure}[H]
    \centering
    \begin{tikzpicture}[xscale=1.5,yscale=1.5]
        \coordinate(O)at(0,0);
        \coordinate(XS)at(-0.5,0);
        \coordinate(XL)at(5.5,0);
        \coordinate(YS)at(0,-0.5);
        \coordinate(YL)at(0,3.7);
        \draw[semithick,->,>=stealth](XS)--(XL)node[right]{$p$};
        \draw[semithick,->,>=stealth](YS)--(YL)node[left]{$q$};
        \coordinate(P)at(1,0);
        \coordinate(Q)at(0,1);

        \foreach\k in{0,1,2,3,4,5}\foreach\l in{0,1,2,3}\fill($(O)+\k*(P)+\l*(Q)$)circle(0.05);
        \foreach\k in{0,1,2,3,4,5}
        \coordinate[label=above right:$1$](X\k)at($\k*(P)$);
        \coordinate[label=below right:$\alpha$](X)at(1,0);
        \foreach\k in{2,3,4,5}
        \coordinate[label=below right:$\alpha^{\k}$](X\k)at($\k*(P)$);
        \coordinate[label=above right:$2g$](Y)at(0,1);
        \foreach\k in{1,2,3,4,5}
        \coordinate[label=above right:$0$](Y\k)at($\k*(P)+(Q)$);
        \coordinate[label=above left:$2g^2+3g+1$, label=below left:$v$](Z)at(0,2);
        \coordinate[label=above:$2g+1$, label=below:$v \alpha$](Z1)at(1,2);
        \foreach\k in{2,3,4,5}\coordinate[label=above right:$1$, label=below right:$v \alpha^{\k}$](Z\k)at($\k*(P)+(Z)$);
        \coordinate(W)at(0,3);
        \foreach\k in{0,1,2,3,4,5}\coordinate[label=above right:$0$](W\k)at($\k*(P)+(W)$);

        \draw(-0.3,3.5)--(0.5,3.5)--(0.5,2.5)--(1.5,2.5)--(1.5,1.5)--(2.5,1.5)--(2.5,0.5)--(3.5,0.5)--(3.5,-0.3);

        \draw [-{Stealth[length=3mm]}] (0.1,1.95) -- (2.9,0.1);
        \coordinate[label=above right:$d_3$](d)at(1.4,1);
        \draw [-{Stealth[length=3mm]}] (2.1,1.95) -- (4.9,0.1);
        \coordinate[label=above right:$d_3$](d)at(3.4,1);

    \end{tikzpicture}
    \caption{$E_3$-term}
    
\end{figure}

To vanish $E_{3}^{2,2}$, the $d_3$ at the right in the figure has to send $v \alpha^{2}$ to $\alpha^{5}$. So the Leibniz rule implies $d_3(v)=\alpha^{3}$. Thus we finally obtain $E_\infty$ as the figure below.
\begin{figure}[H]
    \centering
    \begin{tikzpicture}[xscale=1.5,yscale=1.5]
        \coordinate(O)at(0,0);
        \coordinate(XS)at(-0.5,0);
        \coordinate(XL)at(5.5,0);
        \coordinate(YS)at(0,-0.5);
        \coordinate(YL)at(0,3.7);
        \draw[semithick,->,>=stealth](XS)--(XL)node[right]{$p$};
        \draw[semithick,->,>=stealth](YS)--(YL)node[left]{$q$};
        \coordinate(P)at(1,0);
        \coordinate(Q)at(0,1);

        \foreach\k in{0,1,2,3,4,5}\foreach\l in{0,1,2,3}\fill($(O)+\k*(P)+\l*(Q)$)circle(0.05);
        \foreach\k in{0,1,2}
        \coordinate[label=above right:$1$](X\k)at($\k*(P)$);
        \foreach\k in{3,4,5}
        \coordinate[label=above right:$0$](X\k)at($\k*(P)$);
        \coordinate[label=below right:$\alpha$](X)at(1,0);
        \coordinate[label=below right:$\alpha^{2}$](X2)at($2*(P)$);
        
        \coordinate[label=above right:$2g$](Y)at(0,1);
        \foreach\k in{1,2,3,4,5}
        \coordinate[label=above right:$0$](Y\k)at($\k*(P)+(Q)$);
        \coordinate[label=above left:$2g^2+3g$](Z)at(0,2);
        \coordinate[label=above:$2g$](Z1)at(1,2);
        \foreach\k in{2,3,4,5}\coordinate[label=above right:$0$](Z\k)at($\k*(P)+(Z)$);
        \coordinate(W)at(0,3);
        \foreach\k in{0,1,2,3,4,5}\coordinate[label=above right:$0$](W\k)at($\k*(P)+(W)$);

        \draw(-0.3,3.5)--(0.5,3.5)--(0.5,2.5)--(1.5,2.5)--(1.5,1.5)--(2.5,1.5)--(2.5,0.5)--(3.5,0.5)--(3.5,-0.3);

    \end{tikzpicture}
    \caption{$E_\infty$-term}
    
\end{figure}
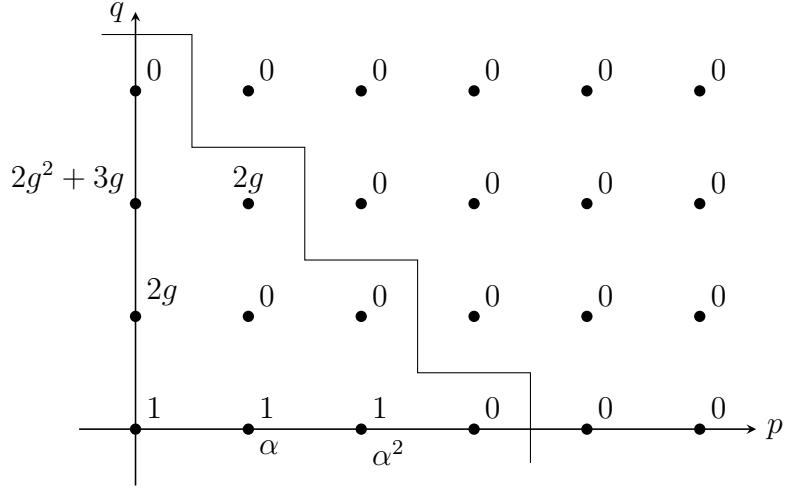

Then we have an isomorphism of $\ccohomology{*}{\rp{\infty}}{\F_2}$-modules
\begin{align}
    \ccohomology{*}{\uconfi{2}{M}}{\F_2}\cong \F_2[\alpha]/(\alpha^3)\oplus \left(\bigoplus_{i=1}^{g}\F_2 x_i\right)\oplus\left(\bigoplus_{i=1}^{2g^2+g}\F_2 y_i \right)\oplus \left(\bigoplus_{i=1}^{2g}\F_2[\alpha]/(\alpha^2)z_i\right).
\end{align}
Here $x_i$ are the basis of $E_{\infty}^{0,1}$, and $y_i$ are the basis of $E_{\infty}^{0,2}$ which derive from $\ccohomology{0}{\sphere{\infty}}{\F_2}\subset E_{2}^{0,2}$. And $z_i$ are the basis of $E_{\infty}^{0,2}$ which derive from $\ccohomology{0}{\rp{\infty}}{\F_2}\subset E_{2}^{0,2}$.

So we have completed the proof of \cref{b_ori}.

\vspace{5mm}
\noindent
\textbf{Case 3}: $M$ is $k$-fold connected sum $\rp{2}\#\cdots \# \rp{2}$.

The cohomology algebra $\ccohomology{*}{M}{\F_2}$ is the quotient of the polynomial ring $\F_2[\alpha_1,\dots,\alpha_k]$ by the relations
\begin{align}
    \alpha_{i}^{3}=0, \quad \alpha_i \alpha_j=0, \quad \alpha_{i}^{2}=\alpha_{j}^{2}=u
\end{align}
where each generator $\alpha_i$ is associated to the dual 1-cell of $i$-th $\rp{2}$. Thus the diagonal class $u_0$ of $M$ is given by $u_0= u\times 1+1\times u+\sum_{i}\alpha_i \times \alpha_i$. So $\ccohomology{2}{\confi{2}{M}}{\F_2}$ is decomposed into
\begin{align}
    \ccohomology{2}{\confi{2}{M}}{\F_2}\cong \ccohomology{2}{M\times M}{\F_2}/\F_2 u_0\cong\F^{\oplus k-1} \oplus \F_2[\sym{2}]^{\oplus \frac{1}{2}k(k+1)+1}
\end{align}
as a $\F_2[\sym{2}]$-module. But the $\F_2$-part is generated by $\alpha_1 \times \alpha_1, \dots ,\alpha_{k-1}\times\alpha_{k-1}$ and the $\F_2[\sym{2}]$-part is generated by $u\times 1$ and $\alpha_i \times \alpha_j (i\neq j)$.

Since $(\alpha_i \times 1)u_0 = u \times \alpha_i + \alpha_i \times u$, we have a direct sum decomposition
\begin{align}
    \ccohomology{3}{\confi{2}{M}}{\F_2}\cong \F_2^{\oplus k}
\end{align}
with generators $u\times \alpha_i(i=1\dots,k)$. The fourth cohomology $\ccohomology{4}{\confi{2}{M}}{\F_2}$ is zero since $\confi{2}{M}$ is open 4-manifold. Thus the proof of \cref{c_unori} has completed.

\vspace{5mm}
By \cref{c_unori}, $E_2$-term is given by
\begin{align}\label{directsumdecomposition2}
    E_{2}^{p,q}&=\ccohomology{p}{\rp{\infty}}{\lcohomology{q}{\confi{2}{M}}{\F_2}}\\
    &=\begin{dcases}
        \ccohomology{p}{\rp{\infty}}{\F_2} & \text{$q=0$,}\\
        \ccohomology{p}{\rp{\infty}}{\F_2[\sym{2}]}^{\oplus k} & \text{$q=1$,}\\
        \ccohomology{p}{\rp{\infty}}{\F_2}^{\oplus k-1}\oplus \ccohomology{p}{\rp{\infty}}{\F_2[\sym{2}]}^{\oplus \frac{1}{2}k(k+1)+1} & \text{$q=2$,}\\
        \ccohomology{p}{\rp{\infty}}{\F_2}^{\oplus k} & \text{$q=3$,}
    \end{dcases}\\
    &=\begin{dcases}
        \ccohomology{p}{\rp{\infty}}{\F_2} & \text{$q=0$,}\\
        \ccohomology{p}{\sphere{\infty}}{\F_2}^{\oplus k} & \text{$q=1$,}\\
        \ccohomology{p}{\rp{\infty}}{\F_2}^{\oplus k-1}\oplus \ccohomology{p}{\sphere{\infty}}{\F_2}^{\oplus \frac{1}{2}k(k+1)+1} & \text{$q=2$,}\\
        \ccohomology{p}{\rp{\infty}}{\F_2}^{\oplus k} & \text{$q=3$.}
    \end{dcases}
\end{align}
So the $E_2$-term is given as the following figure.

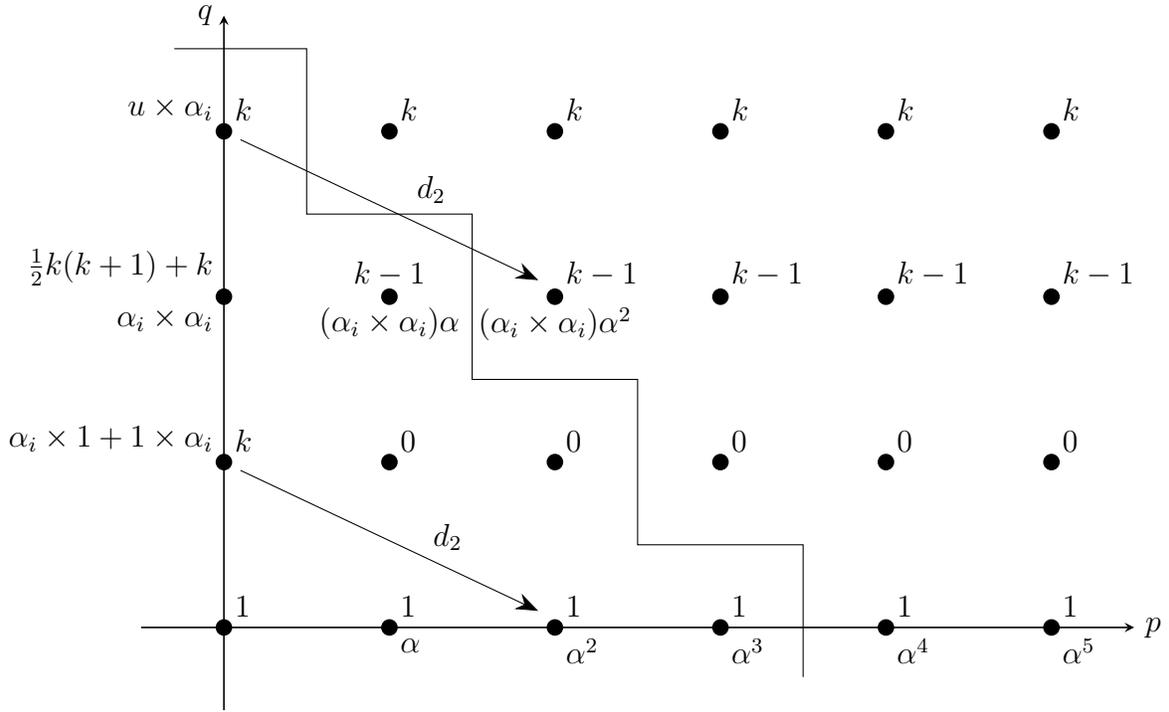
\begin{figure}[H]
    \centering
    \begin{tikzpicture}[xscale=2.2,yscale=2.2]
        \coordinate(O)at(0,0);
        \coordinate(XS)at(-0.5,0);
        \coordinate(XL)at(5.5,0);
        \coordinate(YS)at(0,-0.5);
        \coordinate(YL)at(0,3.7);
        \draw[semithick,->,>=stealth](XS)--(XL)node[right]{$p$};
        \draw[semithick,->,>=stealth](YS)--(YL)node[left]{$q$};
        \coordinate(P)at(1,0);
        \coordinate(Q)at(0,1);

        \foreach\k in{0,1,2,3,4,5}\foreach\l in{0,1,2,3}\fill($(O)+\k*(P)+\l*(Q)$)circle(0.05);
        \foreach\k in{0,1,2,3,4,5}
        \coordinate[label=above right:$1$](X\k)at($\k*(P)$);
        \coordinate[label=below right:$\alpha$](X)at(1,0);
        \foreach\k in{2,3,4,5}
        \coordinate[label=below right:$\alpha^{\k}$](X\k)at($\k*(P)$);
        \coordinate[label=above right:$k$, label=above left:$\alpha_i\times 1+1\times \alpha_i$](Y)at(0,1);
        \foreach\k in{1,2,3,4,5}
        \coordinate[label=above right:$0$](Y\k)at($\k*(P)+(Q)$);
        \coordinate[label=above left:$\frac{1}{2}k(k+1)+k$, label=below left:$\alpha_i\times \alpha_i$](Z)at(0,2);
        \coordinate[label=above:$k-1$, label=below:$(\alpha_i\times \alpha_i)\alpha$](Z1)at(1,2);
        \foreach\k in{2}\coordinate[label=above right:$k-1$, label=below:$(\alpha_i\times\alpha_i) \alpha^{\k}$](Z\k)at($\k*(P)+(Z)$);
        \foreach\k in{3,4,5}\coordinate[label=above right:$k-1$](Z\k)at($\k*(P)+(Z)$);
        \coordinate(W)at(0,3);
        \foreach\k in{0,1,2,3,4,5}\coordinate[label=above right:$k$](W\k)at($\k*(P)+(W)$);
        \coordinate[label=above left:$u\times \alpha_i$](G)at(0,3);

        \draw(-0.3,3.5)--(0.5,3.5)--(0.5,2.5)--(1.5,2.5)--(1.5,1.5)--(2.5,1.5)--(2.5,0.5)--(3.5,0.5)--(3.5,-0.3);

        \draw [-{Stealth[length=3mm]}] (0.1,0.95) -- (1.9,0.1);
        \coordinate[label=above right:$d_2$](d)at(1.2,0.4);
        \draw [-{Stealth[length=3mm]}] (0.1,2.95) -- (1.9,2.1);
        \coordinate[label=above right:$d_2$](d)at(1.1,2.5);

    \end{tikzpicture}
    \caption{$E_2$-term}
    
\end{figure}
The upper $d_2$ is surjective and the lower $d_2$ is zero. So $E_3$-term is described as the figure below.
\begin{figure}[H]
    \centering
    \begin{tikzpicture}[xscale=2.2,yscale=2.2]
        \coordinate(O)at(0,0);
        \coordinate(XS)at(-0.5,0);
        \coordinate(XL)at(5.5,0);
        \coordinate(YS)at(0,-0.5);
        \coordinate(YL)at(0,3.7);
        \draw[semithick,->,>=stealth](XS)--(XL)node[right]{$p$};
        \draw[semithick,->,>=stealth](YS)--(YL)node[left]{$q$};
        \coordinate(P)at(1,0);
        \coordinate(Q)at(0,1);

        \foreach\k in{0,1,2,3,4,5}\foreach\l in{0,1,2,3}\fill($(O)+\k*(P)+\l*(Q)$)circle(0.05);
        \foreach\k in{0,1,2,3,4,5}
        \coordinate[label=above right:$1$](X\k)at($\k*(P)$);
        \coordinate[label=below right:$\alpha$](X)at(1,0);
        \foreach\k in{2,3,4,5}
        \coordinate[label=below right:$\alpha^{\k}$](X\k)at($\k*(P)$);
        \coordinate[label=above right:$k$, label=above left:$\alpha_i\times 1+1\times \alpha_i$](Y)at(0,1);
        \foreach\k in{1,2,3,4,5}
        \coordinate[label=above right:$0$](Y\k)at($\k*(P)+(Q)$);
        \coordinate[label=above left:$\frac{1}{2}k(k+1)+k$, label=below left:$\alpha_i\times \alpha_i$](Z)at(0,2);
        \coordinate[label=above:$k-1$, label=below:$(\alpha_i\times \alpha_i)\alpha$](Z1)at(1,2);
        \foreach\k in{2}\coordinate[label=above right:$0$](Z\k)at($\k*(P)+(Z)$);
        \foreach\k in{3,4,5}\coordinate[label=above right:$0$](Z\k)at($\k*(P)+(Z)$);
        \coordinate(W)at(0,3);
        \foreach\k in{0,1,2,3,4,5}\coordinate[label=above right:$1$](W\k)at($\k*(P)+(W)$);
        \coordinate(G)at(0,3);

        \draw(-0.3,3.5)--(0.5,3.5)--(0.5,2.5)--(1.5,2.5)--(1.5,1.5)--(2.5,1.5)--(2.5,0.5)--(3.5,0.5)--(3.5,-0.3);

        \draw [-{Stealth[length=3mm]}] (0.1,1.95) -- (2.9,0.1);
        \coordinate[label=above right:$d_3$](d)at(1.3,1.1);

    \end{tikzpicture}
    \caption{$E_3$-term}
    
\end{figure}
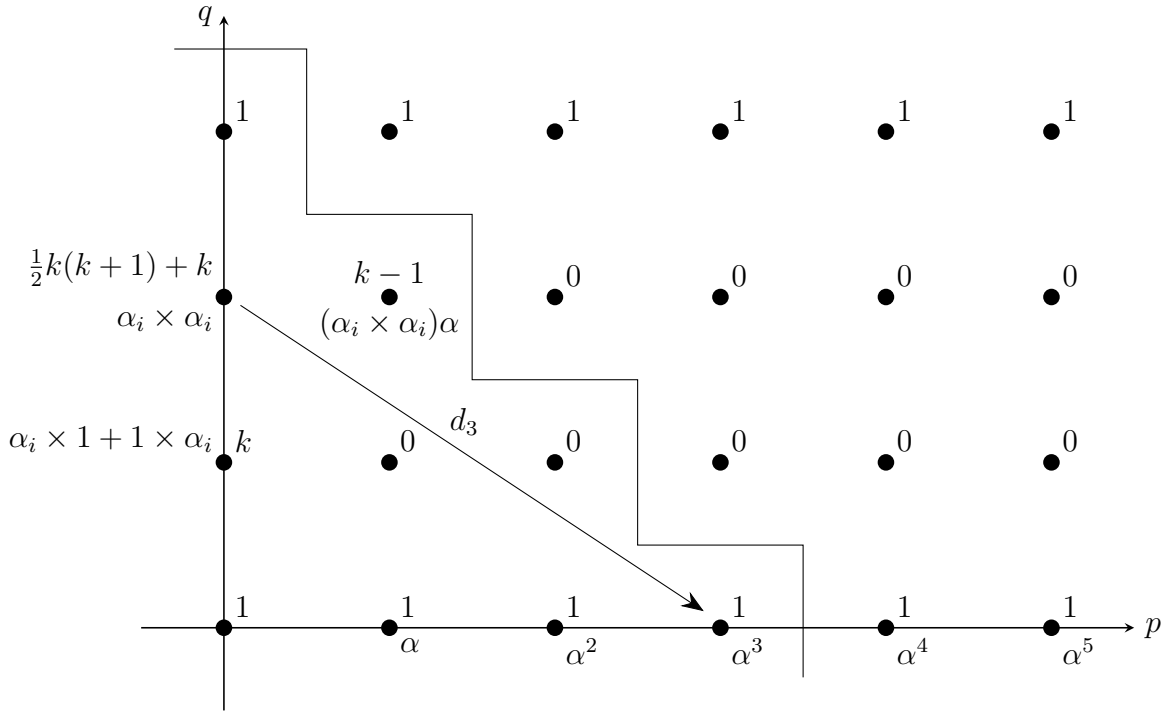
Here we note that $d_3(\alpha_i\times \alpha_i)=0$. To show this, consider the Serre spectral sequence $E^{\prime \prime}$ associated to the Borel construction
\begin{align}\label{borel2}
    \xymatrix{
    M\times M\ar[r]&(M\times M)\times_{\sym{2}}E\sym{2}\ar[r]&B\sym{2}.
    }
\end{align}
Since the fiber of \ref{borel2} has a fixed point $*=(x,x)$, so we obtain an obvious map of fibrations
\begin{align}\label{section}
    \xymatrix{
    M\times M\ar[r]&(M\times M)\times_{\sym{2}}E\sym{2}\ar[r]&B\sym{2}\\
    \ast \ar[r]\ar[u]& \ast \times_{\sym{2}} E\sym{2}=B\sym{2} \ar@{=}[r] \ar[u]& B\sym{2} \ar@{=}[u]
    }
\end{align}
Thus $E_{\infty}^{\prime\prime \ast,0}\cong \ccohomology{*}{\rp{\infty}}{\F_2}=\F_2[\alpha]$. 

$\alpha_i \times \alpha_i \in E_{2}^{\prime \prime0,2}$ is mapped to zero by the derivation since $E_{2}^{\prime \prime 2,1}=0$. So $\alpha_i \times \alpha_i$ survives at $E_{3}^{\prime \prime}$. Further, $\alpha_i \times \alpha_i$ cannot hit to $\alpha^{3}$ at $E_{3}^{\prime \prime}$ because $\alpha^{3}$ survives at $E_{\infty}^{\prime \prime}$. Hence $d_3(\alpha_i\times \alpha_i)=0$ in $E_3$.
\begin{align}
    \xymatrix@C=36pt{
     E_{3}^{\prime\prime}\ar[r]\ar_{d''_3}[d]&E_{3}\ar^{d_3}[d]& \alpha_i\times\alpha_i\ar@{|->}[r]\ar@{|->}[d]&\alpha_i\times\alpha_i\ar@{|->}[d]\\
    E_{3}^{\prime\prime}\ar[r]&E_3 &0\ar@{|->}[r]& 0
    }
\end{align}

Similar to the \textbf{Case 2}, all elements in the $\ccohomology{*}{\sphere{\infty}}{\F_2}$-parts in the decomposition (\ref{directsumdecomposition2}) are mapped to zero by the derivation $d_r$ for any $r\ge 2$. Thus $\alpha^{3}$ survives at $E_4$, and $E_4$ and $E_\infty$ are described as the following. 

\begin{figure}[H]
    \centering
    \begin{tikzpicture}[xscale=1.5,yscale=1.5]
        \coordinate(O)at(0,0);
        \coordinate(XS)at(-0.5,0);
        \coordinate(XL)at(5.5,0);
        \coordinate(YS)at(0,-0.5);
        \coordinate(YL)at(0,3.7);
        \draw[semithick,->,>=stealth](XS)--(XL)node[right]{$p$};
        \draw[semithick,->,>=stealth](YS)--(YL)node[left]{$q$};
        \coordinate(P)at(1,0);
        \coordinate(Q)at(0,1);

        \foreach\k in{0,1,2,3,4,5}\foreach\l in{0,1,2,3}\fill($(O)+\k*(P)+\l*(Q)$)circle(0.05);
        \foreach\k in{0,1,2,3,4,5}
        \coordinate[label=above right:$1$](X\k)at($\k*(P)$);
        \coordinate[label=below right:$\alpha$](X)at(1,0);
        \foreach\k in{2,3,4,5}
        \coordinate[label=below right:$\alpha^{\k}$](X\k)at($\k*(P)$);
        \coordinate[label=above right:$k$](Y)at(0,1);
        \foreach\k in{1,2,3,4,5}
        \coordinate[label=above right:$0$](Y\k)at($\k*(P)+(Q)$);
        \coordinate[label=above left:$\frac{1}{2}k(k+1)+k$](Z)at(0,2);
        \coordinate[label=above:$k-1$](Z1)at(1,2);
        \foreach\k in{2}\coordinate[label=above right:$0$](Z\k)at($\k*(P)+(Z)$);
        \foreach\k in{3,4,5}\coordinate[label=above right:$0$](Z\k)at($\k*(P)+(Z)$);
        \coordinate(W)at(0,3);
        \foreach\k in{0,1,2,3,4,5}\coordinate[label=above right:$1$](W\k)at($\k*(P)+(W)$);
        \coordinate(G)at(0,3);

        \draw(-0.3,3.5)--(0.5,3.5)--(0.5,2.5)--(1.5,2.5)--(1.5,1.5)--(2.5,1.5)--(2.5,0.5)--(3.5,0.5)--(3.5,-0.3);

        \draw [-{Stealth[length=3mm]}] (0.1,2.95) -- (3.9,0.1);
        \coordinate[label=above right:$d_4$](d)at(1.6,1.8);

    \end{tikzpicture}
    \caption{$E_4$-term}

\end{figure}
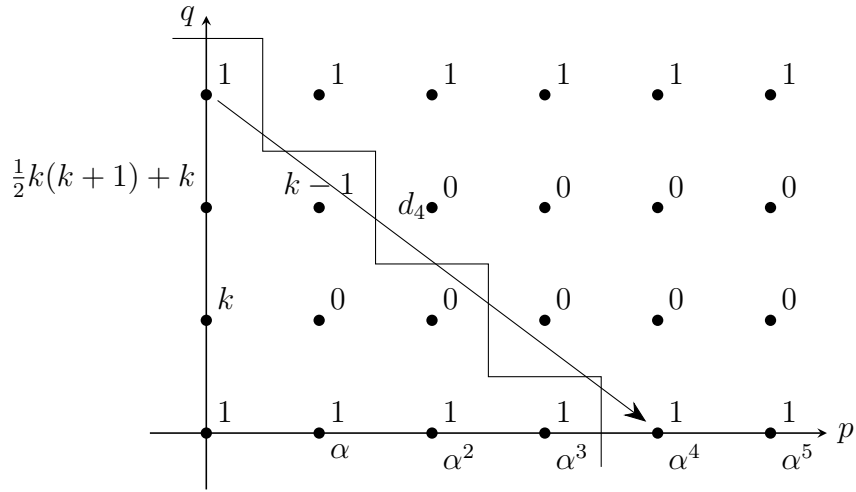
\begin{figure}[H]
    \centering
    \begin{tikzpicture}[xscale=1.5,yscale=1.5]
        \coordinate(O)at(0,0);
        \coordinate(XS)at(-0.5,0);
        \coordinate(XL)at(5.5,0);
        \coordinate(YS)at(0,-0.5);
        \coordinate(YL)at(0,3.7);
        \draw[semithick,->,>=stealth](XS)--(XL)node[right]{$p$};
        \draw[semithick,->,>=stealth](YS)--(YL)node[left]{$q$};
        \coordinate(P)at(1,0);
        \coordinate(Q)at(0,1);

        \foreach\k in{0,1,2,3,4,5}\foreach\l in{0,1,2,3}\fill($(O)+\k*(P)+\l*(Q)$)circle(0.05);
        \foreach\k in{0,1,2,3}
        \coordinate[label=above right:$1$](X\k)at($\k*(P)$);
        \foreach\k in{4,5}
        \coordinate[label=above right:$0$](X\k)at($\k*(P)$);
        \coordinate[label=below right:$\alpha$](X)at(1,0);
        \foreach\k in{2,3}
        \coordinate[label=below right:$\alpha^{\k}$](X\k)at($\k*(P)$);
        \coordinate[label=above right:$k$](Y)at(0,1);
        \foreach\k in{1,2,3,4,5}
        \coordinate[label=above right:$0$](Y\k)at($\k*(P)+(Q)$);
        \coordinate[label=above left:$\frac{1}{2}k(k+1)+k$](Z)at(0,2);
        \coordinate[label=above:$k-1$](Z1)at(1,2);
        \foreach\k in{2}\coordinate[label=above right:$0$](Z\k)at($\k*(P)+(Z)$);
        \foreach\k in{3,4,5}\coordinate[label=above right:$0$](Z\k)at($\k*(P)+(Z)$);
        \coordinate(W)at(0,3);
        \foreach\k in{0,1,2,3,4,5}\coordinate[label=above right:$0$](W\k)at($\k*(P)+(W)$);
        \coordinate(G)at(0,3);

        \draw(-0.3,3.5)--(0.5,3.5)--(0.5,2.5)--(1.5,2.5)--(1.5,1.5)--(2.5,1.5)--(2.5,0.5)--(3.5,0.5)--(3.5,-0.3);

    \end{tikzpicture}
    \caption{$E_\infty$-term}

\end{figure}
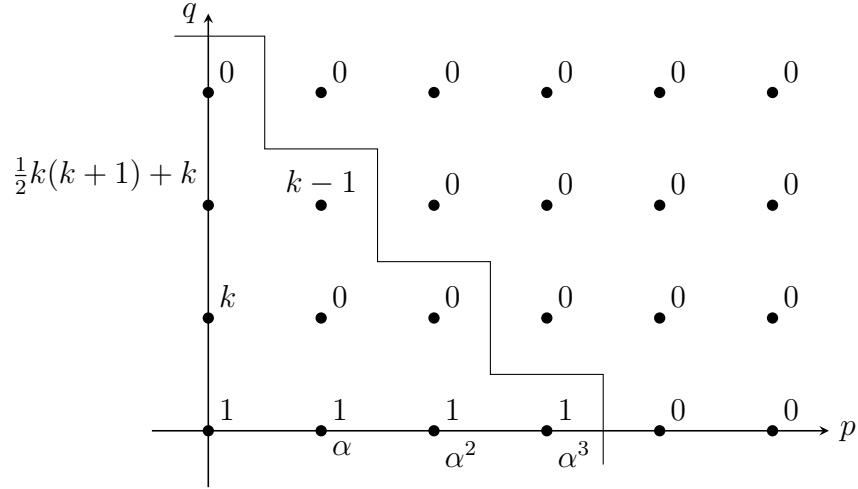

Then we have an isomorphism of $\ccohomology{*}{\rp{\infty}}{\F_2}$-modules
\begin{align}
    \ccohomology{*}{\uconfi{2}{M}}{\F_2}\cong \F_2[\alpha]/(\alpha^4)\oplus \left(\bigoplus_{i=1}^{k}\F_2 x_i\right)\oplus\left(\bigoplus_{i=1}^{\frac{1}{2}k(k+1)+1}\F_2 y_i \right)\oplus \left(\bigoplus_{i=1}^{k-1}\F_2[\alpha]/(\alpha^2)z_i\right).
\end{align}
Here $x_i$ are the basis of $E_{\infty}^{0,1}$, and $y_i$ are the basis of $E_{\infty}^{0,2}$ which derive from $\ccohomology{0}{\sphere{\infty}}{\F_2}\subset E_{2}^{0,2}$. And $z_i$ are the basis of $E_{\infty}^{0,2}$ which derive from $\ccohomology{0}{\rp{\infty}}{\F_2}\subset E_{2}^{0,2}$. Then we have completed the proof of \cref{b_unori}.

\section*{Acknowledgement}
The author appreciates Mitsunobu Tsutaya and Takahiro Matsushita for commenting on this study. 

\section*{Declarations}
The author was also supported by WISE program (MEXT) at Kyushu University. The author has no conflicts of interest relevant to the content of this article.

\bibliographystyle{junsrt}

\end{document}